\documentclass[11pt]{article}
\usepackage{amssymb,amsmath,latexsym, verbatim, enumerate}
\usepackage{tikz}
\usetikzlibrary{shapes}
\allowdisplaybreaks

\begin{document}

\begin{doublespace}

\newtheorem{thm}{Theorem}[section]
\newtheorem{lemma}[thm]{Lemma}
\newtheorem{cond}[thm]{Condition}
\newtheorem{defn}[thm]{Definition}
\newtheorem{prop}[thm]{Proposition}
\newtheorem{corollary}[thm]{Corollary}
\newtheorem{remark}[thm]{Remark}
\newtheorem{example}[thm]{Example}
\newtheorem{conj}[thm]{Conjecture}
\numberwithin{equation}{section}
\def\ee{\varepsilon}
\def\qed{{\hfill $\Box$ \bigskip}}
\def\NN{{\cal N}}
\def\AA{{\cal A}}
\def\MM{{\cal M}}
\def\BB{{\cal B}}
\def\CC{{\cal C}}
\def\LL{{\cal L}}
\def\DD{{\cal D}}
\def\FF{{\cal F}}
\def\EE{{\cal E}}
\def\QQ{{\cal Q}}
\def\RR{{\mathbb R}}
\def\R{{\mathbb R}}
\def\L{{\bf L}}
\def\K{{\bf K}}
\def\S{{\bf S}}
\def\A{{\bf A}}
\def\E{{\mathbb E}}
\def\F{{\bf F}}
\def\P{{\mathbb P}}
\def\N{{\mathbb N}}
\def\eps{\varepsilon}
\def\wh{\widehat}
\def\wt{\widetilde}
\def\pf{\noindent{\bf Proof.} }
\def\beq{\begin{equation}}
\def\eeq{\end{equation}}
\def\lam{\lambda}
\def\H{\mathcal{H}}
\def\nn{\nonumber}
\def\L{\mathcal{L}}

\newcommand\blfootnote[1]{%
  \begingroup
  \renewcommand\thefootnote{}\footnote{#1}%
  \addtocounter{footnote}{-1}%
  \endgroup
}

\title{\Large \bf  Higher order terms of the spectral heat content for killed subordinate and subordinate killed Brownian motions related to symmetric $\alpha$-stable processes in $\R$}
\author{Hyunchul Park}

\date{}
\maketitle


\begin{abstract}
We investigate the 3rd term of the spectral heat content for killed subordinate and subordinate killed Brownian motions on a bounded open interval $D=(a,b)$ in a real line when the underlying subordinators are stable subordinators with index $\alpha \in (1,2)$ or $\alpha=1$. We prove that in the 3rd term of the spectral heat content, one can observe the length $b-a$ of the interval $D$. 
\end{abstract}

\section{Introduction}\label{introduction}
The classical spectral heat content $Q_{D}^{(2)}(t)$ measures the total heat that remains on a domain $D$ with Dirichlet boundary condition and unit initial heat. The spectral heat content can be written in probabilistic terms, and it can be defined as
$$
Q_{D}^{(2)}(t)=\int_{D}\P_{x}(\tau_{D}^{(2)}>t)dx,
$$
where $\tau_{D}^{(2)}=\inf\{t>0 : W_{t}\notin D\}$ is the first exit time from $D$ by a Brownian motion $W=\{W_{t}\}_{t\geq 0}$.
When the Brownian motion is replaced by other L\'evy processes, the corresponding quantity is called the spectral heat content for the L\'evy processes. 
It was recently studied intensively in \cite{V2017, V2016, GPS19}.

One of the most commonly used jump type L\'evy processes is the symmetric stable processes of index $\alpha\in (0,2]$. 
When $\alpha=2$, it is a Brownian motion whose sample paths are continuous with the characteristic exponent being $\E[e^{i\xi W_{t}}]=e^{-t\xi^2}$. 
When $\alpha\in (0,2)$, they are pure-jump processes. 
Stable processes are in fact a special case of subordinate Brownian motions which are time-changed Brownian motions whose time change is given by stable 
subordinators $S_{t}^{(\alpha/2)}$ with Laplace exponent given by
$$
\E[e^{-\lam S_{t}^{(\alpha/2)}}]=e^{-t\lam^{\alpha/2}}, \quad \lam>0.
$$

When one studies the spectral heat content of subordinate Brownian motions, one needs to consider a time-change by a subordinator and killing the process when it first exits the domain under consideration. When we first do time-change and kill the processes, it is called \textit{killed subordinate Brownian motions} and when we first kill the Brownian motions when they first exit the domain and do time-change into the killed Brownian motions, it is called \textit{subordinate killed Brownian motions}. 
These two processes are closely related, and sometimes understanding the spectral heat content of one process helps understand the other. 
The spectral heat content for killed subordinate Brownian motions, when the subordinators are stable subordinators (killed stable processes), were studied in \cite{V2017, V2016}, and the spectral heat content for subordinate killed Brownian motions were studied in \cite{PS19}. In those papers, the authors found the asymptotic expansion of the spectral heat content up to the 2nd terms. 

The purpose of this paper is to refine these results and find the 3rd terms of the spectral heat contents $\tilde{Q}_{D}^{(\alpha)}(t)$ and $Q_{D}^{(\alpha)}(t)$ for subordinate killed Brownian motions and killed subordinate Brownian motions, respectively, in a bounded open interval $D=(a,b)\subset \R$, when the subordinators are stable subordinators for $\alpha\in[1,2)$. 
The main results of this paper are the followings.
The explanation of notations of theorems will be postponed to Section \ref{section:preliminaries} to introduce main results as quickly as possible.
All asymptotic notations are as $t\downarrow 0$.
\begin{thm}\label{thm:KSBM}
Let $D=(a,b)\subset \R$ with $a<b<\infty$, $|D|=b-a$, and $|\partial D|=2$.
\begin{enumerate}[(1)]
\item
Let $\alpha\in (1,2)$. Then, 
\beq\label{eqn:KSBM12}
|D|-Q_{D}^{(\alpha)}(t)
=\E[\overline{X}^{(\alpha)}_{1}]|\partial D|t^{1/\alpha}-\frac{2^{\alpha}\Gamma(\frac{1+\alpha}{2})}{(\alpha-1)\pi^{1/2}\Gamma(1-\frac{\alpha}{2})|D|^{\alpha-1}}t+o(t).
\eeq

\item 
Let $\alpha=1$. Then, 
\begin{align}\label{eqn:KSBM Cauchy}
&|D|-Q_{D}^{(1)}(t) -\frac{1}{\pi}|\partial D|t\ln(\frac{1}{t})\nn\\
=&|\partial D|\left(\int_{0}^{1}\P(\overline{X}^{(1)}_{1}>u)du +\frac{\ln|D|}{\pi} +\int_{1}^{\infty}\P(\overline{X}^{(1)}_{1}>u)-\frac{1}{\pi u}du\right)t +o(t).
\end{align}
\end{enumerate}
\end{thm}

\begin{thm}\label{thm:SKBM}
Let $D=(a,b)\subset \R$ with $a<b<\infty$, $|D|=b-a$, and $|\partial D|=2$.
\begin{enumerate}[(1)]
\item
Let $\alpha\in (1,2)$. Then, 
\beq\label{eqn:SKBM12}
|D|-\tilde{Q}^{(\alpha)}_{D}(t)
=\E[\overline{W}_{S_{1}^{(\frac{\alpha}{2})}}]|\partial D|t^{1/\alpha}-\frac{2\alpha\int_{0}^{\infty}\P(\overline{W}_{1}\geq u)u^{\alpha-1}du}{(\alpha-1)\Gamma(1-\frac{\alpha}{2})|D|^{\alpha-1}} t +o(t).
\eeq

\item 
Let $\alpha=1$. Then, 
\begin{align}\label{eqn:SKBM Cauchy}
&|D|-\tilde{Q}_{D}^{(1)}(t) -\frac{2}{\pi}|\partial D|t\ln(\frac{1}{t})\nn\\
=&|\partial D|\left(\int_{0}^{1}\P(\overline{W}_{S^{(\frac12)}_1}>u)du +\frac{2\ln|D|}{\pi} +\int_{1}^{\infty}\P(\overline{W}_{S^{(\frac12)}_1}>u)-\frac{2}{\pi u} du\right)t+o(t).
\end{align}
\end{enumerate}
\end{thm}
\begin{remark}
When $\alpha\in (0,1)$, the asymptotic expansion for the spectral heat contents $Q_{D}^{(\alpha)}(t)$ or $\tilde{Q}_{D}^{(\alpha)}(t)$ are only known up to the second terms (see \cite[Theorem 1.1]{V2016} and \cite[Theorem 1.1]{PS19}). Asking the third terms when $\alpha\in (0,1)$ is definitely a very interesting question, and we intend to deal with this question in a future project.
\end{remark}

Studying higher order terms is not only an interesting question in itself, but we could also observe that there are some different patterns in the asymptotic expansion of the spectral heat content for Brownian motions and other L\'evy processes by studying higher order terms. 
For Brownian motions, it is well-known that for smooth domains $D$, the spectral heat content has the asymptotic expansion of the form $|D|-Q_{D}^{(2)}(t)\sim \sum_{n=1}^{\infty}a_{n}t^{\frac{n}{2}}$, where $a_{n}$ has some geometric information about the domain $D$ such as perimeter or mean curvature. Hence, it is natural to conjecture that at least when $\alpha\in (1,2)$, the spectral heat content for stable processes is of the form $|D|-Q_{D}^{(\alpha)}(t)\sim \sum_{n=1}^{\infty}b_{n}t^{\frac{n}{\alpha}}$. Theorems \ref{thm:KSBM} and \ref{thm:SKBM} say this is not the case and the asymptotic expansion involves terms that cannot be written as $t^{\frac{n}{\alpha}}$. Also, we observe that the 3rd term involves the length $b-a$ of the underlying interval $D=(a,b)$, hence one can determine the domain $D$ uniquely up to locations, when $D$ is a bounded open interval in $\R$ from the spectral heat content. 

In this paper, we focus on the spectral heat content in dimension one. 
The geometry of open intervals in $\R$ is simple enough to allow to perform detailed computations, and this could be helpful to extend results of this paper into more general settings, such as the spectral heat content in higher dimensions or with respect to more general processes. 
These problems will be studied in forthcoming projects. 

In order to prove the first part of Theorem \ref{thm:KSBM} ($\alpha\in (1,2)$), we analyze the difference $|D|-Q_{D}^{(\alpha)}(t)- \E[\overline{X}^{(\alpha)}_{1}]|\partial D|t^{1/\alpha}$ directly and prove that it is of order $t$. 
Hence, the proof is quite straightforward in this case. 
For the second part of Theorem \ref{thm:KSBM} ($\alpha=1$), the computation becomes delicate because of the logarithmic term $t\ln(1/t)$. 
We utilize the exact form of the density of the supremum process $\overline{X}^{(1)}_{t}=\sup_{s\leq t}X_{s}^{(1)}$ in \cite{Darling56} to compute the difference $\P(\overline{X}^{(1)}_{1}>u)-\frac{1}{\pi u}$ for large $u$, prove that main terms of order $t\ln(1/t)$ cancel out each other, and finally show that the remaining terms are of order $t$.
In order to prove Theorem \ref{thm:SKBM}, we follow a similar path as Theorem \ref{thm:KSBM}. 
For the first part of Theorem \ref{thm:SKBM} ($\alpha\in (1,2)$), we reprove \cite[Theorem 1.1]{PS19} when $D=(a,b)$ and $\alpha\in (1,2)$ using a probabilistic argument in Theorem \ref{thm:prob rep}, which is similar to \cite{V2016}. 
We would like to mention that in Theorem \ref{thm:prob rep}, we express the 2nd coefficient of $|D|-\tilde{Q}_{D}^{(\alpha)}(t)$ by means of the probabilistic term $\E[\overline{W}_{S_{1}^{(\alpha/2)}}]$, which is more natural than other previously known expressions (compare it with \cite[Theorem 1.1]{PS19}). 
In order to prove the second part of Theorem \ref{thm:SKBM} ($\alpha=1$), we establish the tail probability $\P(\overline{W}_{S_{t}^{(\alpha/2)}}>u)$ for $u>1$ in Proposition \ref{prop:tail prob}, which is an amusingly simple expression.
Once having established Proposition \ref{prop:tail prob}, it is straightforward to compute the difference $\P(\overline{W}_{S_{1}^{(\alpha/2)}}>u)-\frac{2}{\pi u}$ for large $u$.
Then, we prove that main terms of order $t\ln(1/t)$ cancel out each other again, and show that the remaining terms are of order $t$. 

The organization of this paper is as follows. 
In Section \ref{section:preliminaries}, we introduce notations and recall some preliminary facts. 
In Section \ref{section:KSBM}, we study the spectral heat content for killed subordinate Brownian motions and prove Theorem \ref{thm:KSBM}. The first part of Theorem \ref{thm:KSBM} is proved in the subsection \ref{subsection:KSBM12}, and the second part of Theorem \ref{thm:KSBM} is proved in the subsection \ref{subsection:KSBM Cauchy}. In Section \ref{section:SKBM}, we study the spectral heat content for subordinate killed Brownian motions, and prove first and second parts of Theorem \ref{thm:SKBM} in subsections \ref{subsection:SKBM12} and \ref{subsection:SKBM Cauchy}, respectively. 
The notation $\P_{x}$ stands for the law of the underlying processes started at $x\in \R$, and $\E_{x}$ stands for expectation with respect to $\P_{x}$. For simplicity, we use $\P=\P_{0}$ and $\E=\E_{0}$.

\section{Preliminaries}\label{section:preliminaries}
In this section, we introduce some notations and define the functions to be studied in the later sections. All stochastic processes and domains are one dimensional objects.

Let $W=\{W_{t}\}_{t\geq 0}$ be a Brownian motion in $\R$. 
The density of the gaussian random variable $W_{t}$ is $p(t,x)=\frac{1}{\sqrt{4\pi t}}e^{-\frac{x^2}{4t}}$ with the characteristic function given by
$$
\E[e^{i\xi W_{t}}]=e^{-t\xi^{2}}, \quad \xi\in\R.
$$
The supremum process $\overline{W}=\{\overline{W}_{t}\}_{t\geq 0}$ of the Brownian motion is defined by $\overline{W}_{t}=\sup_{s\leq t}W_{s}$. 
It follows from \cite[Theorem 2.21]{Le Gall} that $|W_{t}|$ and $\overline{W}_{t}$ have the same distribution.

Let $S^{(\alpha/2)}=\{S_{t}^{(\alpha/2)}\}_{t\geq 0}$ be an $\alpha/2$-stable subordinator. 
That is, $S^{(\alpha/2)}$ is an increasing L\'evy process started at zero whose Laplace exponent is
\beq\label{eqn:SS Laplace exponent}
\E[e^{-\lam S_{t}^{(\alpha/2)}}]=e^{-t\lam^{\alpha/2}}, \quad \lam \geq 0.
\eeq
It follows from \eqref{eqn:SS Laplace exponent} that $S_{t}^{(\alpha/2)}$ and $t^{2/\alpha}S_{1}^{(\alpha/2)}$ have the same distribution for any $t>0$. 
The subordinator $S^{(\alpha/2)}$ is an increasing process started at 0, and for this reason it plays a role as time.
By doing an elementary integral, it is easy to check that
$$
\lam^{\alpha/2} =\frac{\alpha/2}{\Gamma(1-\frac{\alpha}{2})}\int_{0}^{\infty}(1-e^{-\lam t})t^{-1-\frac{\alpha}{2}}dt,
\quad \lam>0, \alpha\in (0,2).
$$
This shows that the L\'evy density $j^{SS}(u)$ for $S^{(\alpha/2)}$ is
\beq\label{eqn:SS Levy density} 
j^{SS}(u)=\frac{\alpha/2}{\Gamma(1-\frac{\alpha}{2})}u^{-1-\frac{\alpha}{2}}, \quad u>0.
\eeq
It follows from \cite[Equation (2.3)]{PS19} or \cite[Equation (18)]{BKKK} that the density $g^{(\alpha/2)}_{1}(x)$ of $S_{1}^{(\alpha/2)}$ exists, and is given by 
\beq\label{eqn:SS density}
g^{(\alpha/2)}_{1}(x)=\frac{1}{\pi}\sum_{n=1}^{\infty}(-1)^{n+1}\frac{\Gamma(1+\frac{\alpha n}{2})}{n!}\sin(\frac{\pi\alpha n}{2})x^{-\frac{\alpha n}{2}-1}, \quad x>0.
\eeq
It follows from the scaling property \eqref{eqn:SS Laplace exponent} that we have 
\beq\label{eqn:SS scaling}
g^{(\alpha/2)}_{t}(x)=t^{-2/\alpha}g^{(\alpha/2)}_{1}(\frac{x}{t^{2/\alpha}}).
\eeq

Now we define subordinate Brownian motions. 
Let $W$ and $S^{(\alpha/2)}$ be Brownian motions and stable subordinators defined on some probability space. 
Assume that they are independent. 
Then, the subordinate Brownian motions by the subordinator $S^{(\alpha/2)}$ are the following time-changed Brownian motions:
$$
X_{t}^{(\alpha)}:=W_{S_{t}^{(\alpha/2)}}.
$$
By conditioning on $S_{t}^{(\alpha/2)}$, one can observe that the characteristic function of time changed process $X^{(\alpha)}=\{X^{(\alpha)}_{t}\}_{t\geq 0}:=\{W_{S_{t}^{(\alpha/2)}}\}_{t\geq 0}$ is given by
\beq\label{eqn:SSP LK}
\E[e^{i\xi X_{t}^{(\alpha)}}]=\E[e^{i\xi W_{S_{t}^{(\alpha/2)}}}]=\E[e^{-S_{t}^{(\alpha/2)}\xi^{2}}]=e^{-t|\xi|^{\alpha}}, \quad \xi \in\R,
\eeq
and this shows that $X^{(\alpha)}$ are symmetric stable processes of index $\alpha$.
From \eqref{eqn:SSP LK}, we observe that $X_{t}^{(\alpha)}$ has the scaling property; $X_{t}^{(\alpha)}$ and $t^{1/\alpha}X_{1}^{(\alpha)}$ have the same distribution for any $t>0$.
The L\'evy density $j^{SSP}(x)$ of $X^{(\alpha)}$ is given by (see \cite[Equation (1.3) and (1.22)]{Bogdan et all})
\beq\label{eqn:SSP Levy density}
j^{SSP}(x)=\frac{A_{1,\alpha}}{|x|^{1+\alpha}}, \quad A_{1,\alpha}=\frac{\alpha2^{\alpha-1}\Gamma(\frac{1+\alpha}{2})}{\pi^{1/2}\Gamma(1-\frac{\alpha}{2})}.
\eeq

Let $D$ be an open set in $\R$, and define $\tau_{D}^{(\alpha)}=\inf\{t>0: X_{t}^{(\alpha)}\notin D\}$ be the first exit time from $D$ by $X^{(\alpha)}$. 
The killed processes $X^{(\alpha), D}=\{X^{(\alpha), D}_{t}\}_{t\geq 0}$ are defined by
$$
X_{t}^{(\alpha), D}=
\begin{cases}
X^{(\alpha)}_{t} & \text{if } t<\tau_{D}^{(\alpha)},\\
\partial & \text{if } t\geq \tau_{D}^{(\alpha)},
\end{cases}
$$
where $\partial$ is a cemetery state. 
The process $X^{(\alpha), D}$ will be called \textit{killed subordinate Brownian motions} (by stable subordinators $S^{(\alpha/2)}$), since we first subordinate (time-change) Brownian motions, then kill the process when they exit the domain. 
We can exchange the order of time-change and killing, and the corresponding process will be called \textit{subordinate killed Brownian motions}  (by stable subordinators $S^{(\alpha/2)}$). 
More precisely, let $\tau_{D}^{(2)}=\inf\{t>0 : W_{t}\notin D\}$ be the first exit time from $D$ by Brownian motions $W$. 
Define killed Brownian motions $W^{D}=\{W^{D}_{t}\}_{t\geq 0}$ as
$$
W_{t}^{D}=
\begin{cases}
W_{t} & \text{if } t<\tau_{D}^{(2)},\\
\partial & \text{if } t\geq \tau_{D}^{(2)}.
\end{cases}
$$
Now the subordinate killed Brownian motions $(W^{D})_{S^{(\alpha/2)}}=\{(W^{D})_{S_{t}^{(\alpha/2)}}\}_{t\geq 0}$ are defined by 
$$
(W^{D})_{S_{t}^{(\alpha/2)}}=
\begin{cases}
W_{S_{t}^{(\alpha/2)}} & \text{if } S_{t}^{(\alpha/2)}<\tau_{D}^{(2)},\\
\partial & \text{if } S_{t}^{(\alpha/2)}\geq \tau_{D}^{(2)}.
\end{cases}
$$
The following graph illustrates sample paths of $X^{(\alpha), D}$ and $(W^{D})_{S^{(\alpha/2)}}$ starting from $x$ when $D=(a,b)$, where the straight lines represent the sample paths of Brownian motions, the circles represent the sample paths of $(W^{D})_{S^{(\alpha/2)}}$, while the circles together with the rectangles represent the sample paths of $X^{(\alpha), D}$.
\begin{figure}
\begin{center}
\begin{tikzpicture}
\draw [->] (-1,0)--(10,0);
\draw [->] (0,-3)--(0,3);
\draw [dotted] (0,2)--(10,2);
\draw [dotted] (0,-2)--(10,-2);
\draw (0,0)--(0.5,1)--(1.5,-1)--(3.5,3)--(6.5,-3)--(8.5,1);
\draw [dotted]  (3,0)--(3,2);
\draw [fill,red] (0,0) circle [radius=0.1];
\draw [fill,red] (0.5,1) circle [radius=0.1];
\draw [fill,red] (1.5,-1) circle [radius=0.1];
\draw [fill,red] (2.75,1.5) circle [radius=0.1];
\filldraw [blue] ([xshift=-2.5pt,yshift=-2.5pt]4.25,1.5) rectangle ++(5pt,5pt);
\filldraw [blue] ([xshift=-2.5pt,yshift=-2.5pt]5,0) rectangle ++(5pt,5pt);
\filldraw [blue] ([xshift=-2.5pt,yshift=-2.5pt]5.75,-1.5) rectangle ++(5pt,5pt);
\node [left] at (0,2) {$b$};
\node [left] at (0,-2) {$a$};
\node [above] at (-0.25,0) {$x$};
\node [below] at (3,0) {$\tau_{D}^{(2)}$};
\end{tikzpicture}
\caption{Sample paths of $(W^{D})_{S^{(\alpha/2)}}$ and $X^{(\alpha), D}$} \label{fig:sample paths}
\end{center}
\end{figure}
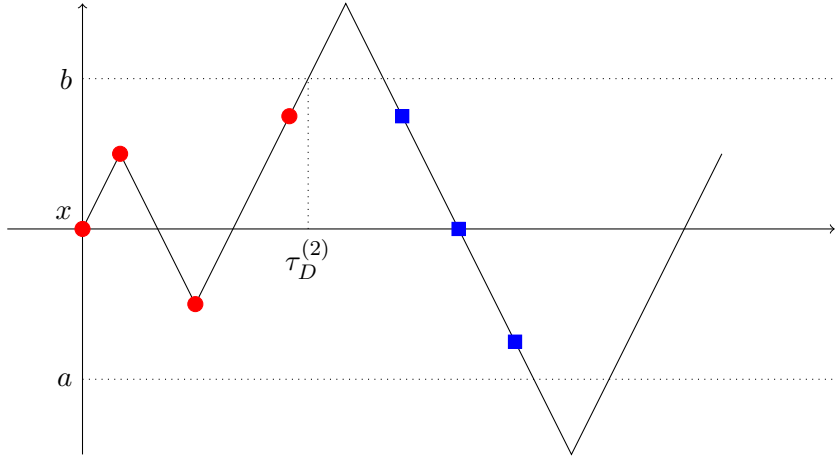
\newpage

Let $\zeta=\inf\{r\geq 0 : (W^{D})_{S_{r}^{(\alpha/2)}}=\partial\}$ be the life time of $(W^{D})_{S_{t}^{(\alpha/2)}}$. 
Then, we have 
$$
\{\zeta >t\} =\{\tau_{D}^{(2)} > S_{t}^{(\alpha/2)} \}.
$$
Clearly, we have $\{\zeta >t\} \subset \{\tau_{D}^{(\alpha)} > t\}$, and the inclusion can be strict. 

We define the supremum processes $\overline{X}^{(\alpha)}=\{\overline{X}_{t}^{(\alpha)}\}_{t\geq 0}$ of the stable processes as  
\beq\label{eqn:sup SS}
\overline{X}^{(\alpha)}_{t}:=\sup_{u\leq t}X^{(\alpha)}_{u}=\sup_{u\leq t}W_{S_{u}^{(\alpha/2)}}.
\eeq
Similarly, $\overline{W}_{S^{(\alpha/2)}}=\{\overline{W}_{S_{t}^{(\alpha/2)}}\}_{t\geq 0}$ are defined by
\beq\label{eqn:sup SB}
\overline{W}_{S^{(\alpha/2)}_{t}}=\sup_{u\leq S^{(\alpha/2)}_{t}}W_{u}.
\eeq
It is noteworthy to mention that even though two expressions $X^{(\alpha)}$ and $W_{S^{(\alpha/2)}}$ mean the same objects, stable processes of index $\alpha$, the supremum notations $\overline{X}^{(\alpha)}$ and $\overline{W}_{S^{(\alpha/2)}}$ are \textit{different}, and we always have $\overline{X}^{(\alpha)}_{t}\leq \overline{W}_{S^{(\alpha/2)}_{t}}$.
The infimum processes $\underline{W}$, $\underline{X}^{(\alpha)}$, and $\underline{W}_{S^{(\alpha/2)}}$ are defined in similar ways with the supremum being replaced by the infimum.

Finally, we define the spectral heat content $Q_{D}^{(\alpha)}(t)$ and $\tilde{Q}_{D}^{(\alpha)}(t)$ for killed subordinate Brownian motions and subordinate killed Brownian motions. 
The spectral heat content $Q_{D}^{(\alpha)}(t)$ for killed subordinate Brownian motions is defined by
$$
Q_{D}^{(\alpha)}(t):=\int_{D}\P_{x}(\tau_{D}^{(\alpha)}>t)dx,
$$
and the spectral heat content $\tilde{Q}_{D}^{(\alpha)}(t)$ for subordinate killed Brownian motions is defined by
$$
\tilde{Q}_{D}^{(\alpha)}(t):=\int_{D}\P_{x}(\zeta>t)dx=\int_{D}\P_{x}(\tau_{D}^{(2)} >S_{t}^{(\alpha/2)})dx.
$$
Since $\{\zeta >t\} \subset \{\tau_{D}^{(\alpha)} > t\}$, we always have 
$$
\tilde{Q}_{D}^{(\alpha)}(t) \leq Q_{D}^{(\alpha)}(t).
$$

When $X^{(\alpha)}$ starts at $x\in (a,b)$, we have 
$$
\{\tau_{D}^{(\alpha)} \leq t\}=\{ \overline{X}^{(\alpha)}_{t} \geq b \text{ or } \underline{X}^{(\alpha)}_{t}\leq a\}.
$$
It follows from the scaling property and the symmetry of $X^{(\alpha)}$, an elementary probability law $\P(A\cup B)=\P(A)+\P(B)-\P(A\cap B)$ for any events $A$ and B, and the change of variable $u=(b-x)t^{-1/\alpha}$ and $v=(a-x)t^{-1/\alpha}$, we have
\begin{align}\label{eqn:1dim expression}
&|D|-Q_{D}^{(\alpha)}(t)=\int_{D}\P_{x}(\tau_{D}^{(\alpha)} \leq t)dx=\int_{a}^{b}\P_{x}(\overline{X}^{(\alpha)}_{t} \geq b \text{ or } \underline{X}^{(\alpha)}_{t}\leq a)dx\nn\\
=&\int_{a}^{b}\P(\overline{X}^{(\alpha)}_{t} \geq b-x)dx+\int_{a}^{b}\P(\underline{X}^{(\alpha)}_{t}\leq a-x)dx -\int_{a}^{b}\P_{x}(\overline{X}^{(\alpha)}_{t} \geq b \text{ and } \underline{X}^{(\alpha)}_{t}\leq a)dx\nn\\
=&\int_{a}^{b}\P(\overline{X}_{1}^{(\alpha)}\geq (b-x)t^{-1/\alpha})dx + \int_{a}^{b}\P(\underline{X}_{1}^{(\alpha)} \leq (a-x)t^{-1/\alpha})dx-\int_{a}^{b}\P_{x}(\overline{X}^{(\alpha)}_{t} \geq b \text{ and } \underline{X}^{(\alpha)}_{t}\leq a)dx\nn\\
=&t^{1/\alpha}\int_{0}^{\frac{b-a}{t^{1/\alpha}}}\P(\overline{X}_{1}\geq u)du +t^{1/\alpha}\int_{-\frac{b-a}{t^{1/\alpha}}}^{0}\P(\underline{X}_{1}^{(\alpha)} \leq v)dv-\int_{a}^{b}\P_{x}(\overline{X}^{(\alpha)}_{t} \geq b \text{ and } \underline{X}^{(\alpha)}_{t}\leq a)dx\nn\\
=&2t^{1/\alpha}\int_{0}^{\frac{b-a}{t^{1/\alpha}}}\P(\overline{X}_{1}\geq u)du-\int_{a}^{b}\P_{x}(\overline{X}^{(\alpha)}_{t} \geq b \text{ and } \underline{X}^{(\alpha)}_{t}\leq a)dx.
\end{align}

\section{Spectral heat content for killed subordinate Brownian motions}\label{section:KSBM}
\subsection{Case: $\alpha\in (1,2)$}\label{subsection:KSBM12}

We start with a simple lemma.
Let $p_{t}^{(\alpha)}(x)$ be the transition density (heat kernel) for $X_{t}^{(\alpha)}$. Note that the following heat kernel estimate is well-known (see \cite{BG});
\beq\label{eqn:HKSS}
c^{-1}(t^{-d/\alpha}\wedge \frac{t}{|x|^{d+\alpha}})\leq p_{t}^{(\alpha)}(x)\leq c(t^{-d/\alpha}\wedge \frac{t}{|x|^{d+\alpha}})
\eeq
for some constant $c>1$.

\begin{lemma}\label{lemma:main term1}
Suppose that $1<\alpha<2$. Then
$$
\lim_{t\rightarrow 0}\frac{t^{1/\alpha}\int_{(b-a)t^{-1/\alpha}}^{\infty}\P(\overline{X}^{(\alpha)}_{1}> u)du}{t}=\frac{2^{\alpha-1}\Gamma(\frac{1+\alpha}{2})}{(\alpha-1)\pi^{1/2}\Gamma(1-\frac{\alpha}{2})(b-a)^{\alpha-1}}.
$$
\end{lemma}
\pf
It follows from L'H\^opital's rule, the scaling property of $X_{t}^{(\alpha)}$, \cite[Proposition VIII 4]{Ber}, \cite[Corollary 8.9]{Sato}, and \eqref{eqn:SSP Levy density} that we have
\begin{align*}
&\lim_{t\rightarrow 0}\frac{\int_{(b-a)t^{-1/\alpha}}^{\infty}\P(\overline{X}^{(\alpha)}_{1}> u)du}{t^{1-1/\alpha}}
=\lim_{t\rightarrow 0}\frac{(b-a)}{\alpha-1}\frac{\P(\overline{X}^{(\alpha)}_{1}>(b-a)t^{-1/\alpha})}{t}\\
=&\lim_{t\rightarrow 0}\frac{(b-a)}{\alpha-1}\frac{\P(X^{(\alpha)}_{1}>(b-a)t^{-1/\alpha})}{t}
=\lim_{t\rightarrow 0}\frac{(b-a)}{\alpha-1}\frac{\P(X^{(\alpha)}_{t}>b-a)}{t}\\
=&\lim_{t\rightarrow 0}\frac{(b-a)}{\alpha-1}\int_{b-a}^{\infty}\frac{p^{(\alpha)}_{t}(u)}{t}du
=\frac{(b-a)}{\alpha-1}\int_{b-a}^{\infty}\frac{A_{1,\alpha}}{u^{1+\alpha}}du=\frac{2^{\alpha-1}\Gamma(\frac{1+\alpha}{2})}{(\alpha-1)\pi^{1/2}\Gamma(1-\frac{\alpha}{2})(b-a)^{\alpha-1}}.
\end{align*}
\qed

\begin{lemma}\label{lemma:remainder1}
Let $\alpha\in (1,2)$. 
Then, for any $t>0$, we have
$$
\int_{a}^{b}\P_{x}(\overline{X}^{(\alpha)}_{t} >b \text{ and } \underline{X}^{(\alpha)}_{t}<a)dx\leq \frac{ct^{1+\frac{1}{\alpha}}}{(b-a)^{\alpha}}\E[\overline{X}^{(\alpha)}_{1}]
$$
for some constant $c>0$.
\end{lemma}
\pf
Define 
$$
\tau:=\inf\{u : X^{(\alpha)}_{u}>b \text{ or } X^{(\alpha)}_{u}<a\}.
$$
Clearly, $\tau$ is a stopping time with respect to the natural filtration $\mathcal{F}=\{\mathcal{F}_{t}\}_{t\geq 0}$.
When the process $X^{(\alpha)}$ starts at $x\in (a,b)$, we have
\begin{align*}
&\{\overline{X}^{(\alpha)}_{t} >b \text{ and } \underline{X}^{(\alpha)}_{t}<a\}
=\{\tau< t, \overline{X}^{(\alpha)}_{t} >b \text{ and } \underline{X}^{(\alpha)}_{t}<a\}\\
=&\{\tau< t, X^{(\alpha)}_{\tau} \leq a,  \overline{X}^{(\alpha)}_{t} >b \text{ and } \underline{X}^{(\alpha)}_{t}<a\}
\cup \{\tau< t, X^{(\alpha)}_{\tau} \geq b,  \overline{X}^{(\alpha)}_{t} >b \text{ and } \underline{X}^{(\alpha)}_{t}<a\}\\
\subset&\{\tau< t, X^{(\alpha)}_{\tau} \leq a,  \overline{X}^{(\alpha)}_{t} >b\}
\cup \{\tau< t, X^{(\alpha)}_{\tau} \geq b,  \underline{X}^{(\alpha)}_{t}<a\}\\
\subset&\{\tau< t, \sup_{\tau \leq s\leq t}(X^{(\alpha)}_{s}-X^{(\alpha)}_{\tau}) >b- a\}
\cup \{\tau< t, \inf_{\tau\leq s \leq t}(X^{(\alpha)}_{s}-X^{(\alpha)}_{\tau}) <-(b-a)\}\\
\subset&\{\tau< t, \sup_{0 \leq s\leq t }(X^{(\alpha)}_{s+\tau}-X^{(\alpha)}_{\tau}) >b- a\}
\cup \{\tau< t, \inf_{0 \leq s \leq t }(X^{(\alpha)}_{s+\tau}-X^{(\alpha)}_{\tau}) <-(b-a)\}\\
=&\{\tau< t, \overline{Y}_{t} >b- a\}\cup \{\tau< t, \underline{Y}_{t} <-(b-a)\},
\end{align*}
where $Y_{u}:=X_{u+\tau}^{(\alpha)}-X_{\tau}^{(\alpha)}$.
By the strong Markov property, $Y$ has the same distribution as $X^{(\alpha)}$ started from 0, and is independent of $\mathcal{F}_{\tau}$.
Hence, for any $x\in (a,b)$, we have
\begin{align}\label{eqn:ub1}
&\P_{x}(\overline{X}^{(\alpha)}_{t} >b \text{ and } \underline{X}^{(\alpha)}_{t}<a)\nn\\
\leq&\P_{x}(\tau< t, \overline{Y}_{t} >b- a ) + \P_{x}(\tau< t, \underline{Y}_{t} <-(b-a))\nn\\
=&\P_{x}(\tau< t)\P(\overline{Y}_{t} >b- a) +\P_{x}(\tau< t)\P( \underline{Y}_{t} <-(b-a))\nn\\
=&\P_{x}(\tau< t)\P(\overline{X}^{(\alpha)}_{t} >b- a) +\P_{x}(\tau< t)\P( \underline{X}^{(\alpha)}_{t} <-(b-a))\nn\\
=&2\P_{x}(\tau< t)\P(\overline{X}^{(\alpha)}_{t} >b- a),
\end{align}
where we used the fact that $\overline{X}^{(\alpha)}_{t}$ and $-\underline{X}^{(\alpha)}_{t}$ have the same distribution because of the symmetry of $X^{(\alpha)}$.

From the scaling property of $X^{(\alpha)}$, \eqref{eqn:HKSS}, and \cite[Proposition 2.1]{V2016}, we have
\beq\label{eqn:ub2}
\P(\overline{X}^{(\alpha)}_{t} > b-a)=\P(\overline{X}^{(\alpha)}_{1}>\frac{b-a}{t^{1/\alpha}})\leq 2\P(X^{(\alpha)}_{1}>\frac{b-a}{t^{1/\alpha}})\leq 
c_1\int_{(b-a)t^{-1/\alpha}}^{\infty}\frac{1}{u^{1+\alpha}}du\leq
\frac{c_{2}t}{(b-a)^{\alpha}}.
\eeq
When $X$ starts at $x\in (a,b)$, we have 
$$
\{\tau<t\}=\{\overline{X}^{(\alpha)}_{t}>b \text{ or } \underline{X}^{(\alpha)}_{t}<a\}.
$$
Hence, from \eqref{eqn:ub1} and \eqref{eqn:ub2}, we have 
\begin{align*}
&\int_{a}^{b}\P_{x}(\overline{X}^{(\alpha)}_{t} >b \text{ and } \underline{X}^{(\alpha)}_{t}<a)dx
\leq\frac{c_{3}t}{(b-a)^{\alpha}}\int_{a}^{b}\P_{x}(\overline{X}^{(\alpha)}_{t}>b \text{ or } \underline{X}^{(\alpha)}_{t}<a)dx\\
\leq&\frac{c_{3}t}{(b-a)^{\alpha}}( \int_{a}^{b}\P_{x}(\overline{X}^{(\alpha)}_{t}>b)dx + \int_{a}^{b}\P_{x}(\underline{X}^{(\alpha)}_{t}<a)dx ).
\end{align*}
By the scaling property of $X$ and the change of variable $u=(b-x)t^{-1/\alpha}$, we have
\begin{align*}
&\int_{a}^{b}\P_{x}(\overline{X}^{(\alpha)}_{t}>b)dx=\int_{a}^{b}\P(\overline{X}^{(\alpha)}_{1}>(b-x)t^{-1/\alpha})
=t^{1/\alpha}\int_{0}^{(b-a)t^{-1/\alpha}}\P(\overline{X}^{(\alpha)}_{1}>u)du
\leq t^{1/\alpha}\E[\overline{X}^{(\alpha)}_{1}].
\end{align*}
Similarly, by the change of variable $v=(x-a)t^{-1/\alpha}$, and the fact that $\overline{X}_{t}$ and $-\underline{X}_{t}$ have the same distribution, we have
\begin{align*}
&\int_{a}^{b}\P_{x}(\underline{X}^{(\alpha)}_{t}<a)dx=\int_{a}^{b}\P(\overline{X}^{(\alpha)}_{1}>(x-a)t^{-1/\alpha})
=t^{1/\alpha}\int_{0}^{(b-a)t^{-1/\alpha}}\P(\overline{X}^{(\alpha)}_{1}>u)du
\leq t^{1/\alpha}\E[\overline{X}^{(\alpha)}_{1}].
\end{align*}
Now the conclusion follows immediately.
\qed

Now we are ready to prove the first part of Theorem \ref{thm:KSBM}.

\noindent{\textbf{Proof of \eqref{eqn:KSBM12}}}

From \eqref{eqn:1dim expression}, we have
\begin{eqnarray}\label{eqn:SHC expansion}
&&|D|-Q_{D}^{(\alpha)}(t)-\E[\overline{X}^{(\alpha)}_{1}]|\partial D|t^{1/\alpha}
=\int_{a}^{b}\P(\tau_{D}^{(\alpha)}\leq t)dx -\E[\overline{X}^{(\alpha)}_{1}]|\partial D|t^{1/\alpha}\nn\\
&=&2t^{1/\alpha}\int_{0}^{\frac{b-a}{t^{1/\alpha}}}\P(\overline{X}^{(\alpha)}_{1}>u)du -\int_{a}^{b}\P_{x}(\overline{X}^{(\alpha)}_{t} >b \text{ and } \underline{X}^{(\alpha)}_{t}<a)dx -2t^{1/\alpha}\int_{0}^{\infty}\P(\overline{X}^{(\alpha)}_{1}>u)du\nn\\
&=&2t^{1/\alpha}\int_{\frac{b-a}{t^{1/\alpha}}}^{\infty}\P(\overline{X}^{(\alpha)}_{1}>u)du -\int_{a}^{b}\P_{x}(\overline{X}^{(\alpha)}_{t} >b \text{ and } \underline{X}_{t}^{(\alpha)}<a)dx.
\end{eqnarray}
Now the conclusion follows immediately from Lemmas \ref{lemma:main term1} and \ref{lemma:remainder1}.  
\qed

\subsection{Case: $\alpha=1$}\label{subsection:KSBM Cauchy}
In this subsection, we study the asymptotic behavior of the spectral heat content for killed subordinate Brownian motions (killed stable processes) when $\alpha=1$. We start with a lemma that is similar to Lemma \ref{lemma:remainder1}.
\begin{lemma}\label{lemma:Cauchy small remainder}
$$
\int_{a}^{b}\P_{x}(\overline{X}^{(1)}_{t}>b \text{ and }\underline{X}^{(1)}_{t}<a )dx=O(t^{2}\ln(1/t)) \text{ as } t\to 0.
$$
\end{lemma}
\pf
The proof is similar to the proof of Lemma \ref{lemma:remainder1}, and we only explain the difference. 
As in the proof of Lemma \ref{lemma:remainder1}, we have 
\begin{eqnarray*}
&&\int_{a}^{b}\P_{x}(\tau<t, \sup_{\tau\leq s \leq t} X^{(1)}_{s}-X^{(1)}_{\tau} >b-a)dx
\leq\frac{ct}{(b-a)^{\alpha}}\int_{a}^{b}\P_{x}(\overline{X}^{(1)}_{t}>b)dx\\
&\leq&\frac{ct^{2}}{(b-a)^{\alpha}}\int_{0}^{\frac{b-a}{t^{1/\alpha}}}\P(\overline{X}^{(1)}_{1}>u)du=O(t^{2}\ln(1/t)),
\end{eqnarray*}
where the last part comes from \cite[Proposition 4.3.(i)]{V2016}.
\qed

There was an error in the paragraph right above \cite[Remark 5.1]{V2016}. The density for $\overline{X}^{(1)}_{1}$ exists and it is given by (see \cite{Darling56})
\beq\label{eqn:sup density}
f(x)=\frac{1}{\pi x^{1/2}(1+x^2)^{3/4}}\exp\left(-\frac{1}{\pi}\int_{0}^{1/x}\frac{\ln v}{1+v^2}dv\right), \quad x>0.
\eeq
We note that there is also a minor error in the exact expression of $f(x)$ in \cite{Darling56} and the upper bound of the integral should be $\frac{1}{x}$, instead of $x$.

Now we are ready to prove the second part of Theorem \ref{thm:KSBM}.

\noindent{\textbf{Proof of \eqref{eqn:KSBM Cauchy}}}

From \eqref{eqn:1dim expression}, we have 
\begin{eqnarray*}
&&|D|-Q^{(1)}_{D}(t)=\int_{a}^{b}\P(\tau^{(1)}_{D}\leq t)dx=2t\int_{0}^{\frac{b-a}{t}}\P(\overline{X}^{(1)}_{1}>u)du  -\int_{a}^{b}\P_{x}(\overline{X}^{(1)}_{t}>b \text{ and }\underline{X}^{(1)}_{t}<a )dx.
\end{eqnarray*}
It follows from Lemma \ref{lemma:Cauchy small remainder}.
$$
\lim_{t\to0}\frac{\int_{a}^{b}\P_{x}(\overline{X}^{(1)}_{t}>b \text{ and }\underline{X}^{(1)}_{t}<a )dx}{t}=0.
$$
Note that from \cite[Proposition 4.3.(i)]{V2016}, we have 
$$
\lim_{t\to 0}\frac{2t\int_{0}^{\frac{b-a}{t}}\P(\overline{X}^{(1)}_{1}>u)du}{t\ln(1/t)}=\frac{2}{\pi}.
$$
We will show that 
\begin{eqnarray*}\label{eqn:Cauchy main 3rd}
&&\lim_{t\to0 }\frac{t\int_{0}^{\frac{b-a}{t}}\P(\overline{X}^{(1)}_{1}>u)du - \frac{t\ln(1/t)}{\pi}}{t}=\lim_{t\to0 }\left(\int_{0}^{\frac{b-a}{t}}\P(\overline{X}^{(1)}_{1}>u)du - \frac{\ln(1/t)}{\pi}\right)\nn\\
&=&\int_{0}^{1}\P(\overline{X}^{(1)}_{1}>u)du +\frac{\ln(b-a)}{\pi} +\int_{1}^{\infty}(\P(\overline{X}^{(1)}_{1}>u)-\frac{1}{\pi u})du.
\end{eqnarray*}

Note for any $0<t<b-a$ that
\begin{eqnarray*}
&&\int_{0}^{\frac{b-a}{t}}\P(\overline{X}^{(1)}_{1}>u)du-\frac{\ln(1/t)}{\pi}\\
&=&\int_{0}^{1}\P(\overline{X}^{(1)}_{1}>u)du+\int_{1}^{\frac{b-a}{t}}\P(\overline{X}^{(1)}_{1}>u)du-\int_{1}^{\frac{b-a}{t}}\frac{1}{\pi u}du +\frac{\ln(b-a)}{\pi}\\
&=&\int_{0}^{1}\P(\overline{X}^{(1)}_{1}>u)du+\frac{\ln(b-a)}{\pi} + \int_{1}^{\frac{b-a}{t}}(\P(\overline{X}^{(1)}_{1}>u) -\frac{1}{\pi u})du.
\end{eqnarray*}
It follows from \eqref{eqn:sup density} and the change of variable $y=\frac{1}{v}$, we have
$$
\P(\overline{X}^{(1)}_{1}>u)=\int_{u}^{\infty}\frac{1}{\pi x^{1/2}(1+x^2)^{3/4}}\exp\left(\frac{1}{\pi}\int_{x}^{\infty}\frac{\ln y }{1+y^2}dy\right)dx.
$$
We will show that for all sufficiently large $u$, we have 
\beq\label{eqn:Cauchy bound}
\left| \P(\overline{X}^{(1)}_{1}>u) -\frac{1}{\pi u}\right|\leq \frac{4}{\pi^2}\frac{\ln u}{u^2},
\eeq
so that by the Lebesgue dominated convergence theorem,
$$
\lim_{t\to 0}\int_{1}^{\frac{b-a}{t}}(\P(\overline{X}^{(1)}_{1}>u) -\frac{1}{\pi u})du=\int_{1}^{\infty}(\P(\overline{X}^{(1)}_{1}>u) -\frac{1}{\pi u})du.
$$

For $u\geq 1$ and $x\geq u$, we have $\exp\left(\frac{1}{\pi}\int_{x}^{\infty}\frac{\ln y }{1+y^2}dy\right)\geq e^{0}=1$ and 
\begin{align*}
&\int_{u}^{\infty}\frac{1}{\pi x^{1/2}(1+x^2)^{3/4}}\exp\left(\frac{1}{\pi}\int_{x}^{\infty}\frac{\ln y }{1+y^2}dy\right)dx
\geq\int_{u}^{\infty}\frac{1}{\pi x^{1/2}(1+x^2)^{3/4}}dx\nn\geq\int_{u}^{\infty}\frac{1}{\pi(1+x^2)}dx\nn\\
=&\frac{1}{\pi}\left(\frac{\pi}{2}-\arctan u\right)=\frac{1}{\pi}\arctan(1/u)\nn=\sum_{n=0}^{\infty}\frac{(-1)^{n}}{\pi (2n+1)}\frac{1}{u^{2n+1}},
\end{align*}
where we used an elementary identity $\arctan u +\arctan(1/u)=\frac{\pi}{2}$.
Hence, there exists $U_{1}>0$ such that 
\beq\label{eqn:Cauchy lb}
\P(\overline{X}^{(1)}_{1}>u)-\frac{1}{\pi u} \geq -\frac{1}{2\pi}\frac{1}{u^3}, \text{ for all } u\geq U_{1}.
\eeq

Now we focus on establishing the upper bound. 
From Karamata's Theorem (\cite[Theorem 1.5.11 (ii)]{BGT}), we have 
$$
\int_{x}^{\infty}\frac{\ln y}{1+y^2}dy=\int_{x}^{\infty}y^{-2}\frac{y^2\ln y}{1+y^2}dy \sim \frac{x\ln x}{1+x^2} \text{ as } x\to \infty.
$$
Hence, there exists $U_{2}>0$ such that for all $x\geq u\geq U_2$, we have 
\beq\label{eqn:Cauchy ub1}
\int_{x}^{\infty}\frac{\ln y}{1+y^2}dy\leq \frac{2x\ln x}{1+x^2}.
\eeq
By an elementary calculus, we see that $e^{u}\leq 1+2u$ for all $0\leq u\leq \ln 2$, and take $U_{3}$ so that 
\beq\label{eqn:Cauchy ub2}
\frac{2}{\pi}\frac{x\ln x}{1+x^2} \leq \ln 2 \text{ for all } x\geq u\geq U_{3}.
\eeq
It follows from \eqref{eqn:Cauchy ub1} and \eqref{eqn:Cauchy ub2} for $u\geq \max(U_2, U_3)$, we have
\begin{eqnarray*}
&&\P(\overline{X}^{(1)}_{1}>u) -\frac{1}{\pi u}\\
&\leq&\int_{u}^{\infty}\frac{1}{\pi x^{1/2}(1+x^2)^{3/4}}\exp\left(\frac{2}{\pi}\frac{x\ln x}{1+x^2}\right)dx -\frac{1}{\pi u}\\
&\leq&\int_{u}^{\infty}\frac{1}{\pi x^{1/2}(1+x^2)^{3/4}}\left( 1+\frac{4}{\pi}\frac{x\ln x}{1+x^2}\right)dx -\frac{1}{\pi u}\\
&\leq&\int_{u}^{\infty}\frac{1}{\pi x^{2}}dx +\int_{u}^{\infty}\frac{1}{\pi x^{1/2}(1+x^2)^{3/4}}\frac{4}{\pi}\frac{x\ln x}{1+x^2}dx-\frac{1}{\pi u}\\
&=&\frac{4}{\pi^{2}}\int_{u}^{\infty}\frac{x^{1/2}\ln x}{(1+x^2)^{7/4}}dx\leq\frac{4}{\pi^{2}}\int_{u}^{\infty}\frac{\ln x}{x^{3}}dx.
\end{eqnarray*}
Again, it follows from \cite[Theorem 1.5.11 (ii)]{BGT}, we have 
$$
\int_{u}^{\infty}\frac{\ln x}{x^{3}}dx \sim \frac12 \frac{\ln u}{u^2} \text{ as } u\to\infty,
$$
and we can take a constant $U_{4} \geq \max(U_2, U_3)$ such that $\int_{u}^{\infty}\frac{\ln x}{x^{3}}dx \leq \frac{\ln u}{u^2}$ for all $u\geq U_{4}$.
Hence, for $u\geq U_{4}$ 
\beq\label{eqn:Cauchy ub3}
\P(\overline{X}^{(1)}_{1}>u) -\frac{1}{\pi u} \leq \frac{4}{\pi^2}\frac{\ln u}{u^2}.
\eeq
Hence, it follows from \eqref{eqn:Cauchy lb} and \eqref{eqn:Cauchy ub3}, there exists $U_{5}\geq \max(U_{1}, U_4)$ such that \eqref{eqn:Cauchy bound} holds for all $u\geq U_{5}$.
\qed

\section{Spectral heat content for subordinate killed Brownian motions}\label{section:SKBM}
In this section, we study the 3rd term of the spectral heat content for subordinate killed Brownian motions, and prove Theorem \ref{thm:SKBM}.

\subsection{Case: $\alpha\in (1,2)$}\label{subsection:SKBM12}

\begin{lemma}\label{lemma:large excursion}
For any $\alpha\in(0,2)$, there exists a constant $c=c(\alpha)>0$ such that 
$$
\P(\overline{W}_{S^{(\alpha/2)}_{t}}>b-a) \leq ct \text{ for all } t>0.
$$
\end{lemma}
\pf
By the scaling property and \cite[Theorem 2.21]{Le Gall}, we have
\begin{eqnarray*}
&&\P(\sup_{u\leq S^{(\alpha/2)}_{t}} W_{u} >b-a)=\P(|W_{S^{(\alpha/2)}_{t}}|>b-a)=\P((S^{(\alpha/2)}_{t})^{1/2}|W_1|>b-a)\\
&=&\P(t^{1/\alpha}(S^{(\alpha/2)}_{1})^{1/2}>\frac{b-a}{|W_1|})=\P(S^{(\alpha/2)}_{1}>\frac{(b-a)^{2}}{t^{2/\alpha}|W_1|^{2}}).
\end{eqnarray*}
Hence, we have 
\begin{eqnarray*}
&&\P(\sup_{u\leq S^{(\alpha/2)}_{t}} W_{u} >b-a)=2\int_{0}^{\infty}\P(S^{(\alpha/2)}_{1}>\frac{(b-a)^2}{t^{2/\alpha}x^2})\frac{1}{\sqrt{4\pi}}e^{-\frac{x^2}{4}}dx\\
&=&2\int_{0}^{\infty}\left(\frac{(b-a)^{2}}{t^{2/\alpha}x^{2}}\right)^{\alpha/2}\P(S^{(\alpha/2)}_{1}>\frac{(b-a)^2}{t^{2/\alpha}x^2})\frac{1}{\sqrt{4\pi}}e^{-\frac{x^2}{4}}\left(\frac{(b-a)^{2}}{t^{2/\alpha}x^{2}}\right)^{-\alpha/2}dx\\
&=&\frac{t}{\sqrt{\pi}(b-a)^{\alpha}}\int_{0}^{\infty}\left(\frac{(b-a)^{2}}{t^{2/\alpha}x^{2}}\right)^{\alpha/2}\P(S^{(\alpha/2)}_{1}>\frac{(b-a)^2}{t^{2/\alpha}x^2}) \times x^{\alpha}e^{-\frac{x^2}{4}}dx.
\end{eqnarray*}
It follows from \cite[Equation (2.8)]{PS19}, there exists a constant $c_1$ such that for all $u\in (0,\infty)$, 
$$
u^{\alpha/2}\P(S^{(\alpha/2)}_{1}>u)\leq c_{1}.
$$
Hence, we have 
$$
\P(\sup_{u\leq S_{t}} W_{u} >b-a)  \leq \frac{c_1 t}{\sqrt{\pi}(b-a)^{\alpha}}\int_{0}^{\infty}x^{\alpha}e^{-\frac{x^2}{4}}dx.
$$
\qed

\begin{lemma}\label{lemma:remainder3}
Let $\alpha\in (1,2)$. Then, there exists a constant $c=c(\alpha)>0$ such that 
$$
\int_{a}^{b}\P_{x}(\sup_{u\leq S_{t}^{(\alpha/2)}}W_{u}>b \text{ and }\inf_{u\leq S_{t}^{(\alpha/2)}}W_{u}<a )dx\leq c\E[\overline{W}_{S_{1}^{(\alpha/2)} }]t^{1+\frac{1}{\alpha}}.
$$
\end{lemma}
\pf
The proof is similar to the proof of Lemma \ref{lemma:remainder1}, and we provide the details for the reader's convenience. 
Define 
$$
\eta:=\inf\{v : \sup_{u\leq S_{v}^{(\alpha/2)}}W_{u}>b \text{ or } \inf_{u\leq S_{v}^{(\alpha/2)}}W_{u}<a\}.
$$
Clearly, $\eta$ is a stopping time with respect to the natural filtration $\mathcal{F}_{t}$.
As in the proof of Lemma \ref{lemma:remainder1}, we have 
\begin{align*}
&\{\sup_{u\leq S_{t}^{(\alpha/2)}}W_{u}>b \text{ and }\inf_{u\leq S_{t}^{(\alpha/2)}}W_{u}<a \}\\
\subset &\{\eta <t, \sup_{u\leq S_{t}^{(\alpha/2)}}\tilde{W}_{u}>b-a\} \cup \{\eta <t,\inf_{u\leq S_{t}^{(\alpha/2)}}\tilde{W}_{u}<-(b-a)  \},
\end{align*}
where $\tilde{W}_{u}:=W_{u+\eta}-W_{\eta}$.
By the strong Markov property, $\tilde{W}$ have the same distribution as $W$ started from 0, and is independent of $\mathcal{F}_{\eta}$. 
Hence, using a similar argument that leads to \eqref{eqn:ub1}, the symmetry of $\tilde{W}$, and Lemma \ref{lemma:large excursion}, we have
\begin{align}\label{eqn:ub new1}
&\P_{x}(\sup_{u\leq S_{t}^{(\alpha/2)}}W_{u}>b \text{ and }\inf_{u\leq S_{t}^{(\alpha/2)}}W_{u}<a )\nn\\
\leq& 2\P_{x}(\eta<t)\P(\sup_{u\leq S_{t}} \tilde{W}_{u} >b-a)
\leq  c_{1}t\P_{x}(\eta<t).
\end{align}
When $W$ starts at $x\in (a,b)$, we have 
$$
\{\eta <t\} = \{  \sup_{u\leq S_{t}^{(\alpha/2)}}W_{u}>b \text{ or } \inf_{u\leq S_{t}^{(\alpha/2)}}W_{u}<a \}.
$$ 
Hence, from \eqref{eqn:ub new1}, we have 
\begin{align}\label{eqn:ub new2}
&\int_{a}^{b}\P_{x}(\sup_{u\leq S_{t}^{(\alpha/2)}}W_{u}>b \text{ and }\inf_{u\leq S_{t}^{(\alpha/2)}}W_{u}<a )dx\nn\\
\leq &c_{1}t\int_{a}^{b}\P_{x}(\sup_{u\leq S_{t}^{(\alpha/2)}}W_{u}>b \text{ or } \inf_{u\leq S_{t}^{(\alpha/2)}}W_{u}<a )dx\nn\\
\leq&c_{1}t(\int_{a}^{b}\P_{x}(\sup_{u\leq S_{t}^{(\alpha/2)}}W_{u}>b)dx + \int_{a}^{b}\P_{x}(\inf_{u\leq S_{t}^{(\alpha/2)}}W_{u}<a )dx).
\end{align}
Note that it follows from the scaling property of $S^{(\alpha/2)}$ and $W$, independence of $S^{(\alpha/2)}$ and $W$, and the change of variable $y=t^{-1/\alpha}(b-x)$, we have
\begin{align}\label{eqn:ub new3}
&\int_{a}^{b}\P_{x}(\sup_{s\leq S^{(\alpha/2)}_{t}}W_{s}\geq b)dx
=\int_{a}^{b}\P(\sup_{s\leq t^{2/\alpha}S^{(\alpha/2)}_{1}}t^{1/\alpha}W_{t^{-2/\alpha}s}\geq b-x)dx\nn\\
=&t^{1/\alpha}\int_{0}^{(b-a)t^{-1/\alpha}}\P(\sup_{u\leq S_{1}^{(\alpha/2)}}W_{u}\geq y)dy
=t^{1/\alpha}\int_{0}^{(b-a)t^{-1/\alpha}}\int_{0}^{\infty}\P(\sup_{u\leq v}W_{u}\geq y)g_{1}^{(\alpha/2)}(v)dvdy\nn\\
\leq&t^{1/\alpha}\int_{0}^{\infty}\int_{0}^{\infty}\P(\sup_{u\leq v}W_{u}\geq y)dy  g_{1}^{(\alpha/2)}(v)dv
= t^{1/\alpha}\int_{0}^{\infty}\E[\overline{W}_{v}]g_{1}^{(\alpha/2)}(v)dv\nn\\
=&t^{1/\alpha}\int_{0}^{\infty}v^{1/2}\E[\overline{W}_1]g_{1}^{(\alpha/2)}(v)dv=t^{1/\alpha}\E[\overline{W}_1]\E[(S_{1}^{(\alpha/2)})^{1/2}]
=t^{1/\alpha}\E[\overline{W}_{S_{1}^{(\alpha/2)} }],
\end{align}
where the last term is known to be finite, since $\alpha>1$ (see \cite[Proposition 2.1]{V2017}).
By the symmetry of $W$, we similarly have 
\begin{align}\label{eqn:ub new4}
&\int_{a}^{b}\P_{x}(\inf_{u\leq S_{t}^{(\alpha/2)}}W_{u}<a )dx
=\int_{a}^{b}\P(\inf_{s\leq t^{2/\alpha}S^{(\alpha/2)}_{1}}t^{1/\alpha}W_{t^{-2/\alpha}s}< a-x)dx\nn\\
=&t^{1/\alpha}\int^{0}_{-(b-a)t^{-1/\alpha}}\P(\inf_{u \leq S_{1}^{(\alpha/2)}}W_{u}<y)dy
=t^{1/\alpha}\int_{0}^{(b-a)t^{-1/\alpha}}\P(\sup_{u\leq S_{1}^{(\alpha/2)}}W_{u}\geq y)dy\nn\\
\leq&t^{1/\alpha}\E[\overline{W}_{S_{1}^{(\alpha/2)} }].
\end{align}
Now the conclusion follows immediately from \eqref{eqn:ub new2}, \eqref{eqn:ub new3}, and \eqref{eqn:ub new4}.
\qed

Now, we reprove the following theorem using the probabilistic argument similar to \cite{V2016}.
\begin{thm}\label{thm:prob rep}
Let $\alpha\in (1,2)$ and $D=(a,b)$ an open interval with $b-a<\infty$. Then, we have
$$
\lim_{t\to 0}\frac{|D|-\tilde{Q}^{(\alpha)}_{D}(t)}{t^{1/\alpha}}=\frac{2\Gamma(1-\frac{1}{\alpha})}{\pi}|\partial D|=\E[\overline{W}_{S^{(\alpha/2)}_{1}}]|\partial D|.
$$
\end{thm}
\pf
The proof is similar to \cite[Theorem 1.1]{V2016}.
When $W$ starts at $x\in (a,b)$, we have
$$
\{\tau_{D}^{(2)}\leq S^{(\alpha/2)}_{t}\}=\{W_{s}\geq b \text{ or } W_{s}\leq a  \text{ for some } s\leq S^{(\alpha/2)}_{t}\}.
$$
Hence, using a similar argument as \eqref{eqn:1dim expression}, we have
\begin{align}\label{eqn:inter1}
&|D|-\tilde{Q}^{(\alpha)}_{D}(t)=\int_{D}\P_{x}(\tau_{D}^{(2)}\leq S_{t}^{(\alpha/2)})dx\nn=\int_{D}\P_{x}(\overline{W}_{S^{(\alpha/2)}_{t}}\geq b \text{ or } \underline{W}_{S^{(\alpha/2)}_{t}}\leq a )dx\nn\\
=&\int_{D}\P_{x}(\overline{W}_{S^{(\alpha/2)}_{t}}\geq b) +\int_{D}\P_{x}(\underline{W}_{S^{(\alpha/2)}_{t}}\leq a)-\int_{D}\P_{x}(\overline{W}_{S^{(\alpha/2)}_{t}}\geq b \text{ and } \underline{W}_{S^{(\alpha/2)}_{t}}\leq a)dx.
\end{align}
From Lemma \ref{lemma:remainder3}, the last expression above is $o(t^{1/\alpha})$.
From \eqref{eqn:ub new3} and \eqref{eqn:ub new4}, and the monotone convergence theorem, we have 
$$
\lim_{t\to 0}\frac{\int_{D}\P_{x}(\overline{W}_{S^{(\alpha/2)}_{t}}\geq b)}{t^{1/\alpha}}
=\lim_{t\to 0}\frac{\int_{D}\P_{x}(\underline{W}_{S^{(\alpha/2)}_{t}}\leq a)}{t^{1/\alpha}}
=\E[\overline{W}_{S_{1}^{(\alpha/2)} }].
$$

Finally, from \cite[Proposition 2.1]{V2017}, we have $\E[(S_{1}^{(\alpha/2)})^{1/2}]=\frac{\Gamma(1-\frac{1}{\alpha})}{\sqrt{\pi}}$ and from \cite[Theorem 2.21]{Le Gall} and a direct computation, we have 
$$
\E[\overline{W}_{1}]=\E[|W_{1}|]=2\int_{0}^{\infty}\frac{x}{\sqrt{4\pi}}e^{-\frac{x^2}{4}}dx=\frac{2}{\sqrt{\pi}}.
$$
From the independence of $W$ and $S^{(\alpha/2)}$, this shows that  
$$
\E[\overline{W}_{S_{1}^{(\alpha/2)} }]=\E[\overline{W}_1]\E[(S_{1}^{(\alpha/2)})^{1/2}]=\frac{2\Gamma(1-\frac{1}{\alpha})}{\pi}.
$$
\qed

Next, we need the following technical computations.
\begin{lemma}\label{lemma:large tail}
We have
$$
\P(\overline{W}_1>u) \sim \frac{2}{\sqrt{\pi}u}e^{-\frac{u^2}{4}} \text{ as } u\to\infty.
$$
\end{lemma}
\pf
It follows from \cite[Theorem 2.21]{Le Gall} that we have
$$
\P(\overline{W}_{1}>u)=\P(|W_1|>u)=2\int_{u}^{\infty}\frac{1}{\sqrt{4\pi}}e^{-\frac{x^2}{4}}dx.
$$
Now it follows from the L'H\^opital's rule, we have 
$$
\lim_{u\to\infty}\frac{2\int_{u}^{\infty}\frac{1}{\sqrt{4\pi}}e^{-\frac{x^2}{4}}d}{ \frac{2}{\sqrt{\pi}u}e^{-\frac{u^2}{4}} }
=\lim_{u\to \infty}\frac{-e^{-\frac{u^2}{4}}}{-\frac{2}{u^2}e^{-\frac{u^2}{4}} -e^{-\frac{u^2}{4}}}=1.
$$
\qed

\begin{lemma}\label{lemma:remainder4}
Let $\alpha\in (1,2)$. Then, we have
$$
\lim_{t\to 0}\frac{\int_{(b-a)t^{-1/\alpha}}^{\infty}\int_{0}^{y^2}\P(\overline{W}_{1}\geq \frac{y}{\sqrt{v}})g_{1}^{(\alpha/2)}(v)dvdy   }{t^{1-\frac{1}{\alpha}}}
=\frac{\alpha}{(\alpha-1)\Gamma(1-\frac{1}{\alpha})(b-a)^{\alpha-1}}\int_{1}^{\infty}\P(\overline{W}_{1}\geq u)u^{\alpha-1}du.
$$
\end{lemma}
\pf
By L'H\^opital's rule, the change of variable $u=\frac{(b-a)t^{-1/\alpha}}{\sqrt{v}}$, \eqref{eqn:SS Levy density}, \eqref{eqn:SS scaling}, \cite[Corollary 8.9]{Sato}, and the Lebesgue dominated convergence theorem using Lemma \ref{lemma:large tail}, we have
\begin{eqnarray*}
&&\lim_{t\to 0}\frac{\int_{(b-a)t^{-1/\alpha}}^{\infty}\int_{0}^{y^2}\P(\overline{W}_{1}\geq \frac{y}{\sqrt{v}})g_{1}^{(\alpha/2)}(v)dvdy   }{t^{1-\frac{1}{\alpha}}}\\
&=&\lim_{t\to 0}\frac{(b-a)}{(\alpha-1)t}\int_{0}^{(b-a)^{2}t^{-2/\alpha}}\P(\overline{W}_{1}\geq \frac{(b-a)t^{-1/\alpha}}{\sqrt{v}})g^{(\alpha/2)}_{1}(v)dv\\
&=&\lim_{t\to 0}\frac{2(b-a)^{3}}{(\alpha-1)t}\int_{1}^{\infty}\P(\overline{W}_{1}\geq u)g^{(\alpha/2)}_{1}(\frac{(b-a)^{2}t^{-2/\alpha}}{u^2})t^{-2/\alpha}u^{-3}du\\
&=&\lim_{t\to 0}\frac{2(b-a)^{3}}{(\alpha-1)}\int_{1}^{\infty}\P(\overline{W}_{1}\geq u)\frac{g^{(\alpha/2)}_{t}(\frac{(b-a)^{2}}{u^2})}{t}u^{-3}du\\
&=&\frac{\alpha}{(\alpha-1)\Gamma(1-\frac{\alpha}{2})(b-a)^{\alpha-1}}\int_{1}^{\infty}\P(\overline{W}_{1}\geq u)u^{\alpha-1}du.
\end{eqnarray*}
\qed

Recall that it follows from \cite[Equation (2.5)]{PS19},
\beq\label{eqn:limit constant}
\lim_{x\to\infty} g^{(\alpha/2)}_{1}(x)x^{1+\frac{\alpha}{2}}=\frac{\alpha}{2\Gamma(1-\frac{\alpha}{2})}.
\eeq
\begin{lemma}\label{lemma:remainder5}
Let $\alpha\in (1,2)$. Then, we have
$$
\lim_{t\to 0}\frac{\int_{(b-a)t^{-1/\alpha}}^{\infty}\int_{y^2}^{\infty}\P(\overline{W}_{1}\geq \frac{y}{v^{1/2}})g_{1}^{(\alpha/2)}(v)dvdy}{t^{1-\frac{1}{\alpha}}}
=\frac{\alpha}{(\alpha-1)\Gamma(1-\frac{\alpha}{2})(b-a)^{\alpha-1}}\int_{0}^{1}\P(\overline{W}_{1}\geq w)w^{\alpha-1}dw.
$$
\end{lemma}
\pf
By the change of variable $w=\frac{y}{v^{1/2}}$, the inner integral in the numerator can be written as 
$$
\int_{y^2}^{\infty}\P(\overline{W}_{1}\geq \frac{y}{v^{1/2}})g_{1}^{(\alpha/2)}(v)dv
=\int_{y^2}^{\infty}\P(\overline{W}_{1}\geq \frac{y}{v^{1/2}})g_{1}^{(\alpha/2)}(v)dv\\
=\int_{0}^{1}\P(\overline{W}_{1}\geq w)g^{(\alpha/2)}_{1}(\frac{y^2}{w^2})\frac{2y^2}{w^3}dw.
$$
Since $g^{(\alpha/2)}_{1}(x)\leq cx^{-1-\frac{\alpha}{2}}$ for $x\geq 1$, the integral is finite.

By the L'H\^opital's rule, the Lebesgue dominated convergence theorem, and \eqref{eqn:limit constant}, we have 
\begin{eqnarray*}
&&\lim_{t\to 0}\frac{\int_{(b-a)t^{-1/\alpha}}^{\infty}\int_{0}^{1}\P(\overline{W}_{1}\geq w)g^{(\alpha/2)}_{1}(\frac{y^2}{w^2})\frac{2y^2}{w^3}dw dy}{t^{1-\frac{1}{\alpha}}}\\
&=&\lim_{t\to 0} \frac{2(b-a)}{(\alpha-1)t}\int_{0}^{1}\P(\overline{W}_{1}\geq w)g^{(\alpha/2)}_{1}( (\frac{(b-a)t^{-1/\alpha}}{w})^{2}) \frac{((b-a)t^{-1/\alpha})^{2}}{w^3}dw \\
&=&\lim_{t\to0}\frac{2}{(\alpha-1)(b-a)^{\alpha-1}}\int_{0}^{1}\P(\overline{W}_{1}\geq w)g^{(\alpha/2)}_{1}((\frac{(b-a)t^{-1/\alpha}}{w})^{2}) (\frac{(b-a)t^{-1/\alpha}}{w})^{2(1+\frac{\alpha}{2})} w^{\alpha-1}dw\\
&=&\frac{2}{(\alpha-1)(b-a)^{\alpha-1}}\int_{0}^{1}\P(\overline{W}_{1}\geq w)\frac{\alpha}{2\Gamma(1-\frac{\alpha}{2})}w^{\alpha-1}dw.
\end{eqnarray*}
\qed

Now we are ready to prove the first part of Theorem \ref{thm:SKBM}.

\noindent{\textbf{Proof of \eqref{eqn:SKBM12}}}

Note that from \eqref{eqn:ub new3}, \eqref{eqn:ub new4}, and \eqref{eqn:inter1}, we have
\begin{align}\label{eqn:3rd term KSBM}
&|D|-\tilde{Q}^{(\alpha)}_{D}(t)\nn\\
=&2t^{1/\alpha}\int_{0}^{\infty}\int_{0}^{\frac{b-a}{t^{1/\alpha}}}\P(\sup_{u\leq v}W_{u}\geq y)dyg_{1}^{(\alpha/2)}(v)dv
-\int_{D}\P_{x}(\overline{W}_{S_{t}^{(\alpha/2)}}>b \text{ and } \underline{W}_{S_{t}^{(\alpha/2)}}<a )dx.
\end{align}
It follows from Lemma \ref{lemma:remainder3}
$$
\int_{D}\P_{x}(\overline{W}_{S_{t}^{(\alpha/2)}}>b \text{ and } \underline{W}_{S_{t}^{(\alpha/2)}}<a)dx=O(t^{1+\frac{1}{\alpha}}).
$$

Now we focus on the first integral in \eqref{eqn:3rd term KSBM}. Note that from the scaling property of $W$, we have
\begin{align*}
&2t^{1/\alpha}\int_{0}^{\infty}\int_{0}^{\frac{b-a}{t^{1/\alpha}}}\P(\sup_{u\leq v}W_{u}\geq y)dyg_{1}^{(\alpha/2)}(v)dv-2t^{1/\alpha}\int_{0}^{\infty}\P(\sup_{u\leq S_{1}^{(\alpha/2)}}W_{u}\geq y)dy\\
=&2t^{1/\alpha}\int_{0}^{\infty}\int_{0}^{\frac{b-a}{t^{1/\alpha}}}\P(\sup_{u\leq v}W_{u}\geq y)dyg_{1}^{(\alpha/2)}(v)dv-2t^{1/\alpha}\int_{0}^{\infty}\int_{0}^{\infty}\P(\sup_{u\leq  v}W_{u}\geq y)g_{1}^{(\alpha/2)}(v)dvdy\\
=&-2t^{1/\alpha}\int_{0}^{\infty}\int_{\frac{b-a}{t^{1/\alpha}}}^{\infty}\P(\overline{W}_{v}\geq y)dyg_{1}^{(\alpha/2)}(v)dv
=-2t^{1/\alpha}\int_{0}^{\infty}\int_{\frac{b-a}{t^{1/\alpha}}}^{\infty}\P(\overline{W}_{1}\geq \frac{y}{\sqrt{v}})dyg_{1}^{(\alpha/2)}(v)dv\\
=&-2t^{1/\alpha}(\int_{\frac{b-a}{t^{1/\alpha}}}^{\infty}\int_{0}^{y^2}\P(\overline{W}_{1}\geq \frac{y}{\sqrt{v}})g_{1}^{(\alpha/2)}(v)dvdy
+\int_{\frac{b-a}{t^{1/\alpha}}}^{\infty}\int_{y^2}^{\infty}\P(\overline{W}_{1}\geq \frac{y}{\sqrt{v}})g_{1}^{(\alpha/2)}(v)dvdy.
\end{align*}
Now the conclusion follows immediately from Lemmas \ref{lemma:remainder4} and \ref{lemma:remainder5}.
\qed

\subsection{Case: $\alpha=1$}\label{subsection:SKBM Cauchy}

In this subsection, we study the spectral heat content for subordinate killed Brownian motions when the underlying subordinator is $S^{(1/2)}$.

\begin{prop}\label{prop:tail prob}
For any $u>1$, we have 
\beq\label{eqn:tail prob}
\P(\overline{W}_{S_{1}^{(1/2)}}>u)=\frac{2}{\pi}\arctan(1/u).
\eeq
\end{prop}
\pf
It follows from \eqref{eqn:SS density} that the density of $S_{1}^{(1/2)}$ is given by
$$
g^{(1/2)}_{1}(x)=\frac{1}{\pi}\sum_{n=1}^{\infty}(-1)^{n+1}\frac{\Gamma(n+\frac12)}{(2n-1)!}x^{-n-\frac12}, \quad x>0.
$$
It is easy to check $(n!)^{2}\leq (2n-1)!$ for all $n\geq 1$ and $\Gamma(n+\frac12)\leq \Gamma(n+1)$ for all $n\geq 2$. Hence, we have 
\begin{align*}
&\left|\sum_{n=1}^{\infty}(-1)^{n+1}\frac{\Gamma(n+\frac12)}{(2n-1)!}x^{-n-\frac12}\right|\leq\sum_{n=1}^{\infty}\frac{\Gamma(n+\frac12)}{(2n-1)!}x^{-n-\frac12}
\leq\Gamma(\frac32)x^{-\frac32} +\sum_{n=2}^{\infty}\frac{\Gamma(n+1)}{(2n-1)!}x^{-n-\frac12}\\
\leq&\Gamma(\frac32)x^{-\frac32} +\sum_{n=2}^{\infty}\frac{x^{-n-\frac12}}{n!}=\Gamma(\frac32)x^{-\frac32} +x^{-\frac12}(e^{1/x}-1-\frac{1}{x})
=O(\frac{1}{x^{3/2}}) \text{ as } x\to \infty.
\end{align*}
Hence, by the Lebesgue dominated convergence theorem, we have
\begin{eqnarray*}
&&\P(S_{1}^{(1/2)}>\frac{u^2}{x^2})=\int_{\frac{u^2}{x^2}}^{\infty}g^{(1/2)}_{1}(v)dv\\
&=&\int_{\frac{u^2}{x^2}}^{\infty}\frac{1}{\pi}\sum_{n=1}^{\infty}(-1)^{n+1}\frac{\Gamma(n+\frac12)}{(2n-1)!}v^{-n-\frac12}dv
=\frac{1}{\pi}\sum_{n=1}^{\infty}(-1)^{n+1}\frac{\Gamma(n+\frac12)}{(2n-1)!}\frac{1}{n-\frac12}\frac{x^{2n-1}}{u^{2n-1}}.
\end{eqnarray*}
It follows from \cite[Theorem 2.21]{Le Gall} that
$$
\P(\overline{W}_{1}>a)=\P(|W_{1}|>a)=2\int_{a}^{\infty}\frac{1}{\sqrt{4\pi}}\int_{a}^{\infty}e^{-\frac{x^2}{4}}dx=\frac{1}{\sqrt{\pi}}\int_{a}^{\infty}e^{-\frac{x^2}{4}}dx.
$$
Hence, we have
\begin{align}\label{eqn:arctan}
&\P(\overline{W}_{S^{(1/2)}_1}>u)=\int_{0}^{\infty}\P(\overline{W}_{y}>u)g_{1}^{(\alpha/2)}(y)dy
=\int_{0}^{\infty}\P(\overline{W}_{1}>\frac{u}{\sqrt{y}})g_{1}^{(\alpha/2)}(y)dy\nn\\
=&\int_{0}^{\infty}\frac{1}{\sqrt{\pi}}\int_{\frac{u}{\sqrt{y}}}^{\infty}e^{-\frac{x^2}{4}}dxg_{1}^{(\alpha/2)}(y)dy
=\frac{1}{\sqrt{\pi}}\int_{0}^{\infty}\int_{\frac{u^2}{x^2}}^{\infty} g_{1}^{(\alpha/2)}(y)dy e^{-\frac{x^2}{4}}dx\nn\\
=&\frac{1}{\sqrt{\pi}}\int_{0}^{\infty}\left(\frac{1}{\pi}\sum_{n=1}^{\infty}(-1)^{n+1}\frac{\Gamma(n+\frac12)}{(2n-1)!}\frac{1}{n-\frac12}\frac{x^{2n-1}}{u^{2n-1}}\right)e^{-\frac{x^2}{4}}dx\nn\\
=&\frac{2}{\pi^{3/2}}\sum_{n=1}^{\infty}(-1)^{n+1}\frac{\Gamma(n+\frac12)}{(2n-1)!}\frac{1}{2n-1}\frac{1}{u^{2n-1}}\int_{0}^{\infty}x^{2n-1}e^{-\frac{x^2}{4}}dx\nn\\
=&\frac{2}{\pi^{3/2}}\sum_{n=1}^{\infty}(-1)^{n+1}\frac{\Gamma(n+\frac12)}{(2n-1)!}\frac{1}{2n-1}\frac{1}{u^{2n-1}}\Gamma(n)2^{2n-1},
\end{align}
where we used $\int_{0}^{\infty}x^{2n-1}e^{-\frac{x^2}{4}}dx=\Gamma(n)2^{2n-1}$, and the interchange of the infinite sum and integral is valid, because of the exponential decay term and the fact $u>1$. 
By the Legendre duplication formula, we have $\Gamma(n)\Gamma(n+\frac12)=2^{1-2n}\sqrt{\pi}\Gamma(2n)$. By the Taylor expansion of $\arctan(x)=\sum_{n=1}^{\infty}(-1)^{n+1}\frac{x^{2n-1}}{2n-1} \text{ for } |x|<1$, \eqref{eqn:arctan} can be simplified to 
$$
\frac{2}{\pi}\sum_{n=1}^{\infty}(-1)^{n+1}\frac{1}{2n-1}\frac{1}{u^{2n-1}}=\frac{2}{\pi}\arctan(\frac{1}{u}),  \quad \text{for }u>1.
$$
\qed
\begin{remark}
Even though it is not necessary for our result, it would be interesting to see if \eqref{eqn:tail prob} holds for all $u>0$. 
\end{remark}

\begin{lemma}\label{lemma:crossing term skbm}
We have
$$
\int_{a}^{b}\P_{x}(\overline{W}_{S^{(1/2)}_t}>b \text{ and } \underline{W}_{S^{(1/2)}_t}<a)dx=O(t^2\ln(1/t)).
$$
\end{lemma}
\pf
The proof is almost identical to the proof of Lemma \ref{lemma:Cauchy small remainder} using Lemma \ref{lemma:large excursion}.
It follows from Proposition \ref{prop:tail prob}, $\P(\overline{W}_{S_{1}^{(1/2)} }>u) \sim \frac{2}{\pi u}$ as $u\to\infty$, and this shows that
$\int_{0}^{\frac{b-a}{t}}\P(\overline{W}_{S_{1}^{(1/2)} }>u)du=O(\ln(1/t))$.
\qed

Now we are ready to prove the second part of Theorem \ref{thm:SKBM}.

\noindent{\textbf{Proof of \eqref{eqn:SKBM Cauchy}}}

Note that from \eqref{eqn:inter1}, we have
\begin{eqnarray*}
&&|D|-\tilde{Q}_{D}^{(1)}(t)=\int_{a}^{b}\P_{x}(\tau_{D}^{(2)}\leq S_{t}^{(1/2)})dx\\
&=&2\int_{a}^{b}\P_{x}(\overline{W}_{S^{(1/2)}_{t}}>b)dx-\int_{a}^{b}\P_{x}(\overline{W}_{S^{(1/2)}_t}>b \text{ and } \underline{W}_{S^{(1/2)}_t}<a)dx.
\end{eqnarray*}
It follows from Lemma \ref{lemma:crossing term skbm}, the second term is $O(t^2\ln(1/t))$.

Now the first expression above can be written as
\begin{eqnarray*}
&&2\int_{a}^{b}\P_{x}(\overline{W}_{S^{(1/2)}_{t}}>b)dx=2\int_{a}^{b}\P_{x}(\sup_{u\leq t^{2}S^{(1/2)}_{1}}tW_{ut^{-2}}>b)dx\\
&=&2\int_{a}^{b}\P_{x}(\sup_{v\leq S^{(1/2)}_{1}}W_{v}>b/t)dx=2t\int_{0}^{\frac{b-a}{t}}\P(\sup_{v\leq S_{1}^{(1/2)}} W_{v}>u )du=2t\int_{0}^{\frac{b-a}{t}}\P(\overline{W}_{S_{1}^{(1/2)}}>u)du.
\end{eqnarray*}
Hence, we have 
\begin{align*}
&2t\int_{0}^{\frac{b-a}{t}}\P(\overline{W}_{S_{1}^{(1/2)}}>u)du-\frac{4}{\pi}t\ln(1/t)=2t(\int_{0}^{\frac{b-a}{t}}\P(\overline{W}_{S^{(1/2)}_1}>u)du -\frac{2}{\pi}\ln(1/t))\\
=&2t(\int_{0}^{1}\P(\overline{W}_{S^{(1/2)}_1}>u)du +\frac{2\ln(b-a)}{\pi} +\int_{1}^{\frac{b-a}{t}}\P(\overline{W}_{S^{(1/2)}_1}>u)-\frac{2}{\pi u} du).
\end{align*}

From Proposition \ref{prop:tail prob}, we have $\P(\overline{W}_{S^{(1/2)}_1}>u)-\frac{2}{\pi u}=O(\frac{1}{u^3})$, and this shows that it is integrable on $(1,\infty)$. Hence, it follows from the monotone convergence theorem 
\begin{eqnarray*}
&&\lim_{t\to 0}\frac{2t}{t}(\int_{0}^{1}\P(\overline{W}_{S^{(1/2)}_1}>u)du +\frac{2\ln(b-a)}{\pi} +\int_{1}^{\frac{b-a}{t}}\P(\overline{W}_{S^{(1/2)}_1}>u)-\frac{2}{\pi u} du)\\
&=&2(\int_{0}^{1}\P(\overline{W}_{S^{(1/2)}_1}>u)du +\frac{2\ln(b-a)}{\pi} +\int_{1}^{\infty}\P(\overline{W}_{S^{(1/2)}_1}>u)-\frac{2}{\pi u} du).
\end{eqnarray*}
\qed

\bigskip
\noindent
{\bf Acknowledgment:}
The author thanks the anonymous referee for carefully reading the manuscript and providing useful suggestions and recommendations.

\bigskip
\noindent

\begin{singlespace}

\end{singlespace}

\end{doublespace}

\vskip 0.3truein

{\bf Hyunchul Park}

Department of Mathematics, State University of New York at New Paltz, NY 12561,
USA

E-mail: \texttt{parkh@newpaltz.edu}

\end{document}